\documentclass [12pt]{article}
\usepackage{graphicx,amssymb,amsfonts,latexsym,amsmath,amsthm,times}
\usepackage{epsfig}
\usepackage{color}
\usepackage[english]{babel}
\usepackage[a4paper]{geometry}
\def\lanbox{\hbox{$\, \vrule height 0.25cm width 0.25cm depth 0.01cm \,$}}

\numberwithin{equation}{section}

\begin{document}

\vspace*{1.4cm}

\normalsize \centerline{\Large \bf ON THE WELL-POSEDNESS
OF SOME MODEL WITH THE }

\medskip

\centerline{\Large \bf CUBED LAPLACIAN ARISING IN THE MATHEMATICAL BIOLOGY}

\vspace*{1cm}

\centerline{\bf Messoud Efendiev$^{1,2}$, Vitali Vougalter$^{3 \ *}$}

\bigskip

\centerline{$^1$ Helmholtz Zentrum M\"unchen, Institut f\"ur Computational
Biology, Ingolst\"adter Landstrasse 1}

\centerline{Neuherberg, 85764, Germany}

\centerline{e-mail: messoud.efendiyev@helmholtz-muenchen.de}

\centerline{$^2$ Azerbaijan University of Architecture and Construction, Baku, Azerbaijan}

\centerline{e-mail: messoud.efendiyev@gmail.com}

\bigskip

\centerline{$^{3 \ *}$  Department of Mathematics, University
of Toronto}

\centerline{Toronto, Ontario, M5S 2E4, Canada}

\centerline{ e-mail: vitali@math.toronto.edu}

\medskip


\vspace*{0.25cm}

\noindent {\bf Abstract:}
In the article we establish the global well-posedness in
$W^{1,(6,2)}({\mathbb R}\times {\mathbb R}^{+})$ of the
integro-differential equation containing the cube of the 
one dimensional Laplacian and the transport
term. Our proof relies on a fixed point technique. 
Furthermore, we formulate the condition leading to the
existence of the nontrivial solution for our problem
under the consideration.
This problem is relevant to the cell
population dynamics in the Mathematical Biology.

\vspace*{0.25cm}



\vspace*{0.5cm}

\bigskip

\bigskip


\setcounter{section}{1}

\centerline{\bf 1. Introduction}

\medskip

\noindent
The present work deals with  the global well-posedness
of the nonlocal reaction-diffusion problem with the constants
$\displaystyle{a\geq 0}$ and $b\in {\mathbb R}$, 
\begin{equation}
\label{h}
\frac{\partial u}{\partial t} =
\frac{\partial^{6}u}{\partial x^{6}}+
b\frac{\partial u}{\partial x}+au+
\int_{-\infty}^{\infty}G(x-y)F(u(y,t), y)dy, \quad x\in {\mathbb R}
\end{equation}
important for the studies of the cell population dynamics. Let us assume that the initial
condition for (\ref{h}) is 
\begin{equation}
\label{ic}
u(x,0)=u_{0}(x)\in H^{6}({\mathbb R}).
\end{equation} 
The analogous equation on ${\mathbb R} $ containing the fractional Laplacian in the context of the anomalous diffusion was studied
in ~\cite{EV25}.
The existence of stationary solutions of the integro-differential problems 
with the bi-Laplacian and the biological applications of such models but without a drift term were covered in ~\cite{VV210}.
The situations on the whole real line and
on a finite interval with periodic boundary conditions involving the drift term and the square root of the one dimensional
negative Laplace operator were discussed in
~\cite{EV22}. The article  ~\cite{EV20} is devoted to the normal diffusion and the transport.
Solvability of certain integro-differential problems with anomalous diffusion, transport
and the cell influx/efflux was established in ~\cite{VV21}.
Spatial structures and generalized travelling waves for an
integro-differential equation were investigated in ~\cite{ABVV10}. Spatial patterns arising in higher order models in physics
and mechanics were studied in ~\cite{PT01}. The work
~\cite{VV130} is about the 
emergence and propagation of patterns in nonlocal reaction-
diffusion equations arising in the theory of speciation and involving the
transport term.  Pattern and waves for a model in population dynamics 
with nonlocal consumptions of resources were considered in ~\cite{GVA06}.
The existence of steady states and travelling waves for the
non-local Fisher-KPP problem was demonstrated in ~\cite{BNPR09}. The work
~\cite{BHN05} is devoted to the estimation of the speed of propagation for KPP type
equations in the periodic framework. Significant applications to the theory of
reaction-diffusion problems involving non-Fredholm operators were developed in
~\cite{DMV05}, ~\cite{DMV08}. Fredholm structures, topological invariants and applications were considered in ~\cite{E09}.
Evolution equations arising in the modelling of life sciences were treated in ~\cite{E13}.
In the article ~\cite{GK18} the authors developed the entropy method for generalized Poisson-Nernst-
Planck equations. The large time behavior of solutions of fourth order parabolic equations and $\epsilon$-entropy of their attractors were  covered in ~\cite{EP07}. Lower estimate of the attractor dimension for a chemotaxis growth system was derived in
~\cite{ATEYM06}. The work ~\cite{E10} is devoted to the development of the theory of finite and infinite dimensional attractors for evolution equations of mathematical physics.
Attractors for degenerate parabolic type problems were studied in ~\cite{E131}.
Exponential decay toward equilibrium via entropy methods 
for reaction-diffusion problems was demonstrated in ~\cite{DF06}.  Local and global existence for nonlocal multispecies 
advection-diffusion models were established in ~\cite{GHLP22}.
Solvability conditions for a linearized Cahn-Hilliard equation of sixth order
were derived in ~\cite{VV12}. Quasilinear elliptic equations on half- and quarter-spaces were covered in ~\cite{DDE13}.
Existence and exact multiplicity for quasilinear elliptic equations in quarter-spaces were discussed in ~\cite{DE17}.

The space variable $x$ in the present work corresponds to the cell genotype,
$u(x,t)$ designates the cell density as a function of the genotype and time.
The right side of (\ref{h}) describes the evolution of the cell density by means of
the cell proliferation, mutations and transport.
The diffusion term is correspondent to the change of genotype
via the small random mutations, and the integral term describes large mutations.
The function $F(u, x)$ denotes the rate of cell birth, which depends on $u$ and $x$
(density dependent proliferation), and the kernel $G(x-y)$ gives
the proportion of newly born cells changing their genotype from $y$ to $x$.
We assume here that it depends on the distance between the genotypes.

The standard Fourier transform in this context is given by
\begin{equation}
\label{ft}
\widehat{\phi}(p):=\frac{1}{\sqrt{2\pi}}\int_{-\infty}^{\infty}\phi(x)e^{-ipx}dx,
\quad p\in {\mathbb R}.
\end{equation}
Evidently, the upper bound
\begin{equation}
\label{fub}  
\|\widehat{\phi}(p)\|_{L^{\infty}({\mathbb R})}\leq \frac{1}{\sqrt{2\pi}}
\|\phi(x)\|_{L^{1}({\mathbb R})}
\end{equation}
holds (see e.g. ~\cite{LL97}). Obviously, (\ref{fub}) yields
\begin{equation}
\label{fub1}  
\|p^{6}\widehat{\phi}(p)\|_{L^{\infty}({\mathbb R})}\leq \frac{1}{\sqrt{2\pi}}
\Big\|\frac{d^{6}\phi}{dx^{6}}\Big\|_{L^{1}({\mathbb R})}.
\end{equation}
We impose the conditions below on the integral kernel involved in 
equation (\ref{h}).

\medskip

\noindent
{\bf Assumption 1.1.}  {\it Suppose $G(x): {\mathbb R}\to {\mathbb R}$ does not vanish identically on the real line.
Moreover,
$\displaystyle{G(x), \ \frac{d^{6}G(x)}{dx^{6}}\in L^{1}({\mathbb R})}$.}     

\medskip

This enables us to introduce the technical quantity
\begin{equation}
\label{g}  
q:=\sqrt{\|G(x)\|_{L^{1}({\mathbb R})}^{2}+
\Big\|\frac{d^{6}G(x)}{dx^{6}}\Big\|_{L^{1}({\mathbb R})}^{2}}.
\end{equation}
Hence, $0<q<\infty$.

\medskip

From the point of view of the applications, the space dimension is
not restricted to $d=1$ because our space variable is correspondent to the cell
genotype but not to the usual physical space.
We use the Sobolev space
\begin{equation}
\label{ss}  
H^{6}({\mathbb R}):=
\Bigg\{\phi(x):{\mathbb R}\to {\mathbb {\mathbb R}} \ | \
\phi(x)\in L^{2}({\mathbb R}), \ \frac{d^{6}\phi}{dx^{6}}\in
L^{2}({\mathbb R}) \Bigg \}.
\end{equation}
It is equipped with the norm
\begin{equation}
\label{n}
\|\phi\|_{H^{6}({\mathbb R})}^{2}:=\|\phi\|_{L^{2}({\mathbb R})}^{2}+
\Bigg\| \frac{d^{6}\phi}{dx^{6}}\Bigg\|_{L^{2}({\mathbb R})}^{2}.
\end{equation}
In order to demonstrate that the global well-posedness of problem (\ref{h}), (\ref{ic}) holds,
we have the function space
$$
W^{1, (6, 2)}({\mathbb R}\times [0, T]):=
$$
\begin{equation}
\label{142}
\Big\{u(x,t): {\mathbb R}\times [0, T]\to {\mathbb R} \ \Big|
\ u(x,t), \ \frac{\partial^{6}u}{\partial x^{6}}, \
\frac{\partial u}{\partial t}
\in L^{2}({\mathbb R}\times [0, T]) \Big\},
\end{equation}  
such that
$$
\|u(x,t)\|_{W^{1, (6, 2)}({\mathbb R}\times [0, T])}^{2}:=
$$
\begin{equation}
\label{142n}
\Big\|\frac{\partial u}{\partial t}\Big\|_{L^{2}({\mathbb R}\times [0, T])}^{2}+
\Big\|\frac{\partial^{6}u}{\partial x^{6}}\Big\|_{L^{2}({\mathbb R}\times [0, T])}^{2}+
\|u\|_{L^{2}({\mathbb R}\times [0, T])}^{2},
\end{equation}  
with $T>0$. In definition (\ref{142n}) we use
$$
\|u\|_{L^{2}({\mathbb R}\times [0, T])}^{2}:=\int_{0}^{T}\int_{-\infty}^{\infty}
|u(x,t)|^{2}dxdt.
$$
Throughout the article we will also have the norm
$$
\|u(x,t)\|_{L^{2}({\mathbb R})}^{2}:=\int_{-\infty}^{\infty}|u(x,t)|^{2}dx.
$$

\medskip

\noindent
{\bf Assumption 1.2.} {\it Function
$F(u, x): {\mathbb R}\times{\mathbb R}\to {\mathbb R}$ satisfies the
Caratheodory condition (see ~\cite{K64}), such that
\begin{equation}
\label{lub}
|F(u, x)|\leq k|u|+h(x) \quad for \quad u\in {\mathbb R}, \quad x\in {\mathbb R}
\end{equation}  
with a constant $k>0$ and
$h(x): {\mathbb R}\to {\mathbb R}^{+}, \ h(x)\in L^{2}({\mathbb R})$.
Furthermore, it is a Lipschitz continuous function, so that
\begin{equation}
\label{lip}
|F(u_{1}, x)-F(u_{2}, x)|\leq l|u_{1}-u_{2}| \quad for \quad any \quad
u_{1, 2}\in {\mathbb R}, \quad x\in {\mathbb R}
\end{equation}
with a constant $l>0$.}

\medskip

Throughout the work ${\mathbb R}^{+}$ denotes the nonnegative semi-axis.
The solvability of a local elliptic equation in a bounded domain in
${\mathbb R}^{N}$ was discussed in ~\cite{BO86}. The nonlinear function contained there
was allowed to have a sublinear growth.

Let us apply the standard Fourier transform (\ref{ft}) to both sides of
problem (\ref{h}), (\ref{ic}). This yields
\begin{equation}
\label{hf}
\frac{\partial \widehat{u}}{\partial t}=[-p^{6}+ibp+a]\widehat{u}+
\sqrt{2\pi}\widehat{G}(p)\widehat{f}_{u}(p,t),
\end{equation}
\begin{equation}
\label{icf}
\widehat{u(x,0)}(p)=\widehat{u_{0}}(p).  
\end{equation}
In formula (\ref{hf}) and further down $\widehat{f}_{u}(p,t)$ will stand for the Fourier image
of $F(u(x,t), x)$. Note that
$$
u(x,t)=\frac{1}{\sqrt{2\pi}}\int_{-\infty}^{\infty}\widehat{u}(p,t)e^{ipx}dp,
\quad 
\frac{\partial u}{\partial t}=\frac{1}{\sqrt{2\pi}}\int_{-\infty}^{\infty}
\frac{\partial \widehat{u}(p,t)}{\partial t}e^{ipx}dp
$$
with $x\in {\mathbb R}, \ t\geq 0$. The Duhamel's principle enables us
to reformulate problem (\ref{hf}), (\ref{icf}) as
$$
\widehat{u}(p,t)=
$$
\begin{equation}
\label{duh}
e^{t\{-p^{6}+ibp+a\}}\widehat{u_{0}}(p)+\int_{0}^{t}
e^{(t-s)\{-p^{6}+ibp+a\}}\sqrt{2\pi}\widehat{G}(p)\widehat{f}_{u}(p,s)ds.
\end{equation}
Let us write down the auxiliary equation related to (\ref{duh}), namely
$$
\widehat{u}(p,t)=
$$
\begin{equation}
\label{aux}
e^{t\{-p^{6}+ibp+a\}}\widehat{u_{0}}(p)+\int_{0}^{t}
e^{(t-s)\{-p^{6}+ibp+a\}}\sqrt{2\pi}\widehat{G}(p)\widehat{f}_{v}(p,s)ds,
\end{equation}
where $\widehat{f}_{v}(p,s)$ denotes the Fourier image of
$F(v(x,s), x)$ under transform (\ref{ft}) and $v(x,t)\in W^{1, (6, 2)}({\mathbb R}\times [0, T])$ is arbitrary.

We introduce the operator $t_{a, b}$, such that $u =t_{a, b}v$, where $u$ 
is a solution of  (\ref{aux}). The main proposition of the work is as follows.

\bigskip

\noindent
{\bf Theorem 1.3.} {\it Let Assumptions 1.1 and 1.2 hold and
\begin{equation}
\label{qlt}
ql\sqrt{T^{2}e^{2aT}(1+2[a+|b|+1]^{2})+2}<1, 
\end{equation}
where the constant $q$ is defined in (\ref{g}) and the Lipschitz constant
$l$ is introduced in (\ref{lip}).
Then problem (\ref{aux}) defines the map
$t_{a, b}: W^{1, (6, 2)}({\mathbb R}\times [0, T])\to W^{1, (6, 2)}
({\mathbb R}\times [0, T])$, which is a strict contraction.
The unique fixed point $w(x,t)$ of this map $t_{a, b}$ is the only solution of
(\ref{h}), (\ref{ic}) in
$W^{1, (6, 2)}({\mathbb R}\times [0, T])$.}

\medskip

The final statement of the article deals with  the global well-posedness for
our problem. 

\medskip

\noindent
{\bf Corollary 1.4.} {\it Let the conditions of Theorem 1.3 be fulfilled.
Then (\ref{h}), (\ref{ic}) admits a unique solution   
$w(x,t)\in W^{1, (6, 2)}({\mathbb R}\times {\mathbb R}^{+})$.  This solution does not vanish identically
 for $x\in {\mathbb R}$ and $t\in {\mathbb R}^{+}$ 
provided the intersection of supports of the Fourier transforms of functions
$\hbox{supp}\widehat{F(0, x)}\cap \hbox{supp}\widehat{G}$ is a set of
nonzero Lebesgue measure on ${\mathbb R}$.}

\medskip

We turn our attention to the proof of the main result of the work.

\bigskip


\setcounter{section}{2}
\setcounter{equation}{0}

\centerline{\bf 2. The well-posedness of the equation}

\bigskip

\noindent
{\it Proof of Theorem 1.3.} Let us choose arbitrarily
$v(x,t)\in W^{1, (6, 2)}({\mathbb R}\times [0, T])$ and verify
that the first term in the right side of (\ref{aux})
belongs to $L^{2}({\mathbb R}\times [0, T])$. Clearly,
$$
\|e^{t\{-p^{6}+ibp+a\}}\widehat{u_{0}}(p)\|_{L^{2}({\mathbb R})}^{2}=
\int_{-\infty}^{\infty}e^{-2tp^{6}}e^{2at}|\widehat{u_{0}}(p)|^{2}dp\leq
e^{2at}\|u_{0}\|_{L^{2}({\mathbb R})}^{2}.
$$
Therefore,
$$
\|e^{t\{-p^{6}+ibp+a\}}\widehat{u_{0}}(p)\|_{L^{2}({\mathbb R}\times [0, T])}^{2}=
\int_{0}^{T}\|e^{t\{-p^{6}+ibp+a\}}\widehat{u_{0}}(p)\|_{L^{2}({\mathbb R})}^{2}dt
\leq
$$
$$
\int_{0}^{T}e^{2at}\|u_{0}\|_{L^{2}({\mathbb R})}^{2}dt.
$$
The right side of the last inequality equals to
$\displaystyle{\frac{e^{2aT}-1}{2a}\|u_{0}\|_{L^{2}({\mathbb R})}^{2}}$
when $a>0$ and
$T\|u_{0}\|_{L^{2}({\mathbb R})}^{2}$ if $a=0$. This means that
\begin{equation}
\label{u0l2}  
e^{t\{-p^{6}+ibp+a\}}\widehat{u_{0}}(p)\in L^{2}({\mathbb R}\times [0, T]).
\end{equation}
Obviously, we have the estimate from above on the norm of the second term in the right side of (\ref{aux}) as
$$
\Big\|\int_{0}^{t}
e^{(t-s)\{-p^{6}+ibp+a\}}\sqrt{2\pi}\widehat{G}(p)\widehat{f}_{v}(p,s)ds\Big\|_
{L^{2}({\mathbb R})}\leq
$$
$$  
\int_{0}^{t}\Big\|
e^{(t-s)\{-p^{6}+ibp+a\}}\sqrt{2\pi}\widehat{G}(p)\widehat{f}_{v}(p,s)\Big\|_
{L^{2}({\mathbb R})}ds.
$$
Note that 
$$  
\Big\|
e^{(t-s)\{-p^{6}+ibp+a\}}\sqrt{2\pi}\widehat{G}(p)\widehat{f}_{v}(p,s)\Big\|_
{L^{2}({\mathbb R})}^{2}=
$$
\begin{equation}
\label{intgfv}  
\int_{-\infty}^{\infty}e^{-2(t-s)p^{6}}e^{2a(t-s)}2\pi|\widehat{G}(p)|^{2}
|\widehat{f}_{v}(p,s)|^{2}dp.
\end{equation}
By means of inequality (\ref{fub}) we obtain the upper bound on the right side of
(\ref{intgfv}) as
$$
e^{2a(t-s)}2\pi\|\widehat{G}(p)\|_{L^{\infty}({\mathbb R})}^{2}
\|F(v(x,s),x)\|_{L^{2}({\mathbb R})}^{2}\leq
$$
$$
e^{2aT}\|G(x)\|_{L^{1}({\mathbb R})}^{2}
\|F(v(x,s),x)\|_{L^{2}({\mathbb R})}^{2}.
$$
Thus,
$$  
\Big\|
e^{(t-s)\{-p^{6}+ibp+a\}}\sqrt{2\pi}\widehat{G}(p)\widehat{f}_{v}(p,s)\Big\|_
{L^{2}({\mathbb R})}\leq
$$
$$
e^{aT}\|G(x)\|_{L^{1}({\mathbb R})}
\|F(v(x,s),x)\|_{L^{2}({\mathbb R})}.
$$
Let us recall (\ref{lub}). We derive
\begin{equation}
\label{lubs}  
\|F(v(x,s),x)\|_{L^{2}({\mathbb R})}\leq k\|v(x,s)\|_{L^{2}({\mathbb R})}+\|h(x)\|_
{L^{2}({\mathbb R})}.
\end{equation}
This yields
$$  
\Big\|
e^{(t-s)\{-p^{6}+ibp+a\}}\sqrt{2\pi}\widehat{G}(p)\widehat{f}_{v}(p,s)\Big\|_
{L^{2}({\mathbb R})}\leq
$$
$$
e^{aT}\|G(x)\|_{L^{1}({\mathbb R})}
\{k\|v(x,s)\|_{L^{2}({\mathbb R})}+\|h(x)\|_{L^{2}({\mathbb R})}\},
$$
such that
$$  
\Big\|\int_{0}^{t}
e^{(t-s)\{-p^{6}+ibp+a\}}\sqrt{2\pi}\widehat{G}(p)\widehat{f}_{v}(p,s)ds\Big\|_
{L^{2}({\mathbb R})}\leq
$$
$$
ke^{aT}\|G(x)\|_{L^{1}({\mathbb R})}\int_{0}^{T}\|v(x,s)\|_{L^{2}({\mathbb R})}ds+
Te^{aT}\|G(x)\|_{L^{1}({\mathbb R})}\|h(x)\|_{L^{2}({\mathbb R})}.
$$
Using the Schwarz inequality, we arrive at
\begin{equation}
\label{sch}  
\int_{0}^{T}\|v(x,s)\|_{L^{2}({\mathbb R})}ds\leq
\sqrt{\int_{0}^{T}\|v(x,s)\|_{L^{2}({\mathbb R})}^{2}ds}\sqrt{T}.
\end{equation}
This gives us
$$
\Big\|\int_{0}^{t}
e^{(t-s)\{-p^{6}+ibp+a\}}\sqrt{2\pi}\widehat{G}(p)\widehat{f}_{v}(p,s)ds\Big\|_
{L^{2}({\mathbb R})}^{2}\leq
$$
$$
e^{2aT}\|G(x)\|_{L^{1}({\mathbb R})}^{2}
\{k\sqrt{T}\|v(x,s)\|_{L^{2}({\mathbb R}\times [0,T])}+T\|h(x)\|_{L^{2}({\mathbb R})}\}^{2}.
$$
Let us obtain the estimate from above on the norm as
$$
\Big\|\int_{0}^{t}
e^{(t-s)\{-p^{6}+ibp+a\}}\sqrt{2\pi}\widehat{G}(p)\widehat{f}_{v}(p,s)ds\Big\|_
{L^{2}({\mathbb R}\times [0,T])}^{2}=
$$
$$
\int_{0}^{T}\Big\|\int_{0}^{t}
e^{(t-s)\{-p^{6}+ibp+a\}}\sqrt{2\pi}\widehat{G}(p)\widehat{f}_{v}(p,s)ds\Big\|_
{L^{2}({\mathbb R})}^{2}dt\leq 
$$
$$
e^{2aT}\|G(x)\|_{L^{1}({\mathbb R})}^{2}
\{k\|v(x,s)\|_{L^{2}({\mathbb R}\times [0,T])}+\sqrt{T}\|h(x)\|_{L^{2}({\mathbb R})}\}^{2}
T^{2},
$$
which is finite under the given conditions for
$v(x,s)\in W^{1, (6, 2)}({\mathbb R}\times [0, T])$. Hence,
\begin{equation}
\label{0tetsgvl2}  
\int_{0}^{t}
e^{(t-s)\{-p^{6}+ibp+a\}}\sqrt{2\pi}\widehat{G}(p)\widehat{f}_{v}(p,s)ds\in
L^{2}({\mathbb R}\times [0,T]).
\end{equation}
By means of (\ref{u0l2}), (\ref{0tetsgvl2}) and (\ref{aux}), we derive
\begin{equation}
\label{upt12}  
\widehat{u}(p,t)\in L^{2}({\mathbb R}\times [0,T]),
\end{equation}
so that
\begin{equation}
\label{uxtl2}  
u(x,t)\in L^{2}({\mathbb R}\times [0,T]).
\end{equation}
According to (\ref{aux}),
$$
p^{6}\widehat{u}(p,t)=
$$
\begin{equation}
\label{aux2}
e^{t\{-p^{6}+ibp+a\}}p^{6}\widehat{u_{0}}(p)+\int_{0}^{t}
e^{(t-s)\{-p^{6}+ibp+a\}}\sqrt{2\pi}p^{6}\widehat{G}(p)\widehat{f}_{v}(p,s)ds.
\end{equation}
Let us analyze the first term in the right side of (\ref{aux2}). Evidently,
$$
\|e^{t\{-p^{6}+ibp+a\}}p^{6}\widehat{u_{0}}(p)\|_{L^{2}({\mathbb R}\times [0,T])}^{2}=
\int_{0}^{T}\int_{-\infty}^{\infty}e^{-2tp^{6}}e^{2at}|p^{6}\widehat{u_{0}}(p)|^{2}dp
dt\leq
$$
$$
\int_{0}^{T}\int_{-\infty}^{\infty}e^{2at}|p^{6}\widehat{u_{0}}(p)|^{2}dpdt.
$$
Clearly, this is equal to
$\displaystyle
{\frac{e^{2aT}-1}{2a}\Big\|\frac{d^{6}u_{0}}{dx^{6}}\Big\|_{L^{2}({\mathbb R})}^{2}}$ if
$a>0$ and
$\displaystyle{T\Big\|\frac{d^{6}u_{0}}{dx^{6}}\Big\|_{L^{2}({\mathbb R})}^{2}}$
for $a=0$. Thus,
\begin{equation}
\label{p2u0hpl2}
e^{t\{-p^{6}+ibp+a\}}p^{6}\widehat{u_{0}}(p)\in L^{2}({\mathbb R}\times [0,T]).
\end{equation}  
Consider the second term in the right side of
(\ref{aux2}). Obviously,
$$
\Big\|\int_{0}^{t}
e^{(t-s)\{-p^{6}+ibp+a\}}\sqrt{2\pi}p^{6}\widehat{G}(p)\widehat{f}_{v}(p,s)ds
\Big\|_{L^{2}({\mathbb R})}\leq 
$$
$$
\int_{0}^{t}\Big\|
e^{(t-s)\{-p^{6}+ibp+a\}}\sqrt{2\pi}p^{6}\widehat{G}(p)\widehat{f}_{v}(p,s)
\Big\|_{L^{2}({\mathbb R})}ds.
$$
Evidently,
$$
\Big\|e^{(t-s)\{-p^{6}+ibp+a\}}\sqrt{2\pi}p^{6}\widehat{G}(p)\widehat{f}_{v}(p,s)
\Big\|_{L^{2}({\mathbb R})}^{2}=
$$
\begin{equation}
\label{etsp2gpl2}
\int_{-\infty}^{\infty}
e^{-2(t-s)p^{6}}e^{2a(t-s)}2\pi|p^{6}\widehat{G}(p)|^{2}|\widehat{f}_{v}(p,s)|^{2}
dp.
\end{equation}  
Let us use (\ref{fub1}) to obtain the upper bound on 
the right side of equality (\ref{etsp2gpl2}) as
$$
2\pi e^{2aT}\|p^{6}\widehat{G}(p)\|_{L^{\infty}({\mathbb R})}^{2}
\int_{-\infty}^{\infty}|\widehat{f}_{v}(p,s)|^{2}dp\leq
$$
$$
e^{2aT}\Big\|\frac{d^{6}G}{dx^{6}}\Big\|_{L^{1}({\mathbb R})}^{2}
\|F(v(x,s),x)\|_{L^{2}({\mathbb R})}^{2}.
$$
We recall (\ref{lubs}) . Hence,
$$
\Big\|e^{(t-s)\{-p^{6}+ibp+a\}}\sqrt{2\pi}p^{6}\widehat{G}(p)\widehat{f}_{v}(p,s)
\Big\|_{L^{2}({\mathbb R})}\leq
$$
$$
e^{aT}\Big\|\frac{d^{6}G}{dx^{6}}\Big\|_{L^{1}({\mathbb R})}\{k\|v(x,s)\|_
{L^{2}({\mathbb R})}+\|h(x)\|_{L^{2}({\mathbb R})}\},
$$
such that
$$
\Big\|\int_{0}^{t}e^{(t-s)\{-p^{6}+ibp+a\}}\sqrt{2\pi}p^{6}\widehat{G}(p)
\widehat{f}_{v}(p,s)ds\Big\|_{L^{2}({\mathbb R})}\leq
$$
$$
ke^{aT}\Big\|\frac{d^{6}G}{dx^{6}}\Big\|_{L^{1}({\mathbb R})}\int_{0}^{T}\|v(x,s)\|_
{L^{2}({\mathbb R})}ds+
Te^{aT}\Big\|\frac{d^{6}G}{dx^{6}}\Big\|_{L^{1}({\mathbb R})}\|h(x)\|_{L^{2}({\mathbb R})}.
$$
By virtue of inequality (\ref{sch}), we have
$$
\Big\|\int_{0}^{t}e^{(t-s)\{-p^{6}+ibp+a\}}\sqrt{2\pi}p^{6}\widehat{G}(p)
\widehat{f}_{v}(p,s)ds\Big\|_{L^{2}({\mathbb R})}^{2}\leq
$$
$$
e^{2aT}\Big\|\frac{d^{6}G}{dx^{6}}\Big\|_{L^{1}({\mathbb R})}^{2}
\{k\|v(x,s)\|_{L^{2}({\mathbb R}\times [0,T])}\sqrt{T}+
\|h(x)\|_{L^{2}({\mathbb R})}T\}^{2}.
$$
Therefore,
$$
\Big\|\int_{0}^{t}e^{(t-s)\{-p^{6}+ibp+a\}}\sqrt{2\pi}p^{6}\widehat{G}(p)
\widehat{f}_{v}(p,s)ds\Big\|_{L^{2}({\mathbb R}\times [0, T])}^{2}\leq
$$
$$
e^{2aT}\Big\|\frac{d^{6}G}{dx^{6}}\Big\|_{L^{1}({\mathbb R})}^{2}
\{k\|v(x,s)\|_{L^{2}({\mathbb R}\times [0,T])}+
\|h(x)\|_{L^{2}({\mathbb R})}\sqrt{T}\}^{2}T^{2}<\infty
$$
under the stated assumptions  with
$v(x,s)\in W^{1, (6, 2)}({\mathbb R}\times [0, T])$. This means that
\begin{equation}
\label{int0tp2Gpfv}  
\int_{0}^{t}e^{(t-s)\{-p^{6}+ibp+a\}}\sqrt{2\pi}p^{6}\widehat{G}(p)
\widehat{f}_{v}(p,s)ds\in L^{2}({\mathbb R}\times [0, T]).
\end{equation}
(\ref{p2u0hpl2}), (\ref{int0tp2Gpfv}) and (\ref{aux2}) imply that
\begin{equation}
\label{up2t12}  
p^{6}\widehat{u}(p,t)\in L^{2}({\mathbb R}\times [0, T]).
\end{equation}  
Thus,
\begin{equation}
\label{d2udx2l2}
\frac{\partial^{6}u}{\partial x^{6}}\in L^{2}({\mathbb R}\times [0, T]).
\end{equation}  
Let us recall (\ref{aux}). We derive
\begin{equation}
\label{duhdt}
\frac{\partial \widehat{u}}{\partial t}=\{-p^{6}+ibp+a\}\widehat{u}(p,t)+
\sqrt{2\pi}\widehat{G}(p)\widehat{f}_{v}(p,t).
\end{equation}
By means of (\ref{upt12}),
\begin{equation}
\label{aupt12}  
a\widehat{u}(p,t)\in L^{2}({\mathbb R}\times [0,T]).
\end{equation}
We arrive at the estimate from above on the norm using (\ref{upt12}) and (\ref{up2t12}) as
$$
\|ibp\widehat{u}(p,t)\|_{L^{2}({\mathbb R}\times [0,T])}^{2}=b^{2}\int_{0}^{T}
\Big\{\int_{|p|\leq 1}p^{2}|\widehat{u}(p,t)|^{2}dp+
\int_{|p|>1}p^{2}|\widehat{u}(p,t)|^{2}dp\Big\}dt\leq
$$
$$
b^{2}\{\|\widehat{u}(p,t)\|_{L^{2}({\mathbb R}\times [0,T])}^{2}+
\|p^{6}\widehat{u}(p,t)\|_{L^{2}({\mathbb R}\times [0,T])}^{2}\},
$$
which is finite, such that
\begin{equation}
\label{bupt12} 
ibp\widehat{u}(p,t)\in L^{2}({\mathbb R}\times [0,T]).
\end{equation}
Clearly,
\begin{equation}
\label{uptalf12} 
-p^{6}\widehat{u}(p,t)\in L^{2}({\mathbb R}\times [0,T])
\end{equation}
due to  (\ref{up2t12}).
Combining  (\ref{aupt12}), (\ref{bupt12}) and (\ref{uptalf12}) gives us
\begin{equation}
\label{uptalfab12}
(-p^{6}+ibp+a)\widehat{u}(p,t)\in L^{2}({\mathbb R}\times [0,T]).
\end{equation}
Let us consider the last term in the right side of (\ref{duhdt}). By virtue of inequalities (\ref{fub}) and 
(\ref{lubs}), we have
$$
\|\sqrt{2\pi}\widehat{G}(p)\widehat{f}_{v}(p,t)\|_
{L^{2}({\mathbb R}\times [0,T])}^{2}\leq 2\pi
\|\widehat{G}(p)\|_{L^{\infty}({\mathbb R})}^{2}\int_{0}^{T}
\|F(v(x,t), x)\|_{L^{2}({\mathbb R})}^{2}dt\leq
$$
$$
\|G(x)\|_{L^{1}({\mathbb R})}^{2}\int_{0}^{T}(k\|v(x,t)\|_{L^{2}({\mathbb R})}+
\|h(x)\|_{L^{2}({\mathbb R})})^{2}dt\leq
$$
$$
\|G(x)\|_{L^{1}({\mathbb R})}^{2}\{2k^{2}\|v(x,t)\|_{L^{2}({\mathbb R}\times [0, T])}^{2}+
2\|h(x)\|_{L^{2}({\mathbb R})}^{2}T\}<\infty
$$
under the given conditions for
$v(x,t)\in W^{1, (6, 2)}({\mathbb R}\times [0, T])$. Therefore,
\begin{equation}
\label{ghfvhpt}  
\sqrt{2\pi}\widehat{G}(p)\widehat{f}_{v}(p,t)\in L^{2}({\mathbb R}\times [0,T]).
\end{equation}
By means of equation (\ref{duhdt}) along with (\ref{uptalfab12}) and
(\ref{ghfvhpt}), we obtain
$$
\frac{\partial \widehat{u}}{\partial t}\in L^{2}({\mathbb R}\times [0,T]),
$$
which yields
\begin{equation}
\label{dudtl2}  
\frac{\partial u}{\partial t}\in L^{2}({\mathbb R}\times [0,T]).
\end{equation}
Let us recall the definition of the norm (\ref{142n}). By virtue of the statements
(\ref{uxtl2}), (\ref{d2udx2l2}) and (\ref{dudtl2}), 
$$
u(x,t)\in W^{1, (6, 2)}({\mathbb R}\times [0, T]).
$$
Such function is uniquely determined by (\ref{aux}).
Hence, under the stated assumptions equation (\ref{aux}) defines a map
$$
t_{a, b}: W^{1, (6, 2)}({\mathbb R}\times [0, T])\to
W^{1, (6, 2)}({\mathbb R}\times [0, T]).
$$
Our goal is to demonstrate that under the given conditions this map is a
strict contraction. Let us choose arbitrarily
$v_{1, 2}(x,t)\in W^{1, (6, 2)}({\mathbb R}\times [0, T])$. By means of the
reasoning above,
$u_{1, 2}:=t_{a, b}v_{1, 2}\in W^{1, (6, 2)}({\mathbb R}\times [0, T])$.
Formula  (\ref{aux}) yields
$$
\widehat{u_{1}}(p,t)=
$$
\begin{equation}
\label{1aux}
e^{t\{-p^{6}+ibp+a\}}\widehat{u_{0}}(p)+\int_{0}^{t}
e^{(t-s)\{-p^{6}+ibp+a\}}\sqrt{2\pi}\widehat{G}(p)\widehat{f}_{v_{1}}(p,s)ds,
\end{equation}
$$
\widehat{u_{2}}(p,t)=
$$
\begin{equation}
\label{2aux}
e^{t\{-p^{6}+ibp+a\}}\widehat{u_{0}}(p)+\int_{0}^{t}
e^{(t-s)\{-p^{6}+ibp+a\}}\sqrt{2\pi}\widehat{G}(p)\widehat{f}_{v_{2}}(p,s)ds.
\end{equation}
Here $\widehat{f}_{v_{j}}(p,s), \ j=1, 2$ stands for the Fourier image of
$F(v_{j}(x,s), x)$ under transform (\ref{ft}). According to the system of equations
(\ref{1aux}), (\ref{2aux}),
$$
\widehat{u_{1}}(p,t)-\widehat{u_{2}}(p,t)=
$$
\begin{equation}
\label{u1u2hint0t}  
\int_{0}^{t}
e^{(t-s)\{-p^{6}+ibp+a\}}\sqrt{2\pi}\widehat{G}(p)
[\widehat{f}_{v_{1}}(p,s)-\widehat{f}_{v_{2}}(p,s)]ds.
\end{equation}
Evidently, the estimate from above on the norm 
$$
\|\widehat{u_{1}}(p,t)-\widehat{u_{2}}(p,t)\|_{L^{2}({\mathbb R})}\leq
$$
\begin{equation}
\label{int0tetsf12}  
\int_{0}^{t}\|
e^{(t-s)\{-p^{6}+ibp+a\}}\sqrt{2\pi}\widehat{G}(p)
[\widehat{f}_{v_{1}}(p,s)-\widehat{f}_{v_{2}}(p,s)]\|_{L^{2}({\mathbb R})}ds
\end{equation}
is valid.
By virtue of (\ref{fub}), we have the upper bound 
$$
\|e^{(t-s)\{-p^{6}+ibp+a\}}\sqrt{2\pi}\widehat{G}(p)
[\widehat{f}_{v_{1}}(p,s)-\widehat{f}_{v_{2}}(p,s)]\|_{L^{2}({\mathbb R})}^{2}=
$$
$$
2\pi \int_{-\infty}^{\infty}e^{-2(t-s)p^{6}}e^{2(t-s)a}|\widehat{G}(p)|^{2}
|\widehat{f}_{v_{1}}(p,s)-\widehat{f}_{v_{2}}(p,s)|^{2}dp\leq 
$$
$$
2\pi e^{2aT}\|\widehat{G}(p)\|_{L^{\infty}({\mathbb R})}^{2}\int_{-\infty}^{\infty}
|\widehat{f}_{v_{1}}(p,s)-\widehat{f}_{v_{2}}(p,s)|^{2}dp\leq 
$$
$$
e^{2aT}\|G(x)\|_{L^{1}({\mathbb R})}^{2}\|F(v_{1}(x,s), x)-F(v_{2}(x,s), x)\|_
{L^{2}({\mathbb R})}^{2}.
$$
Let us recall condition (\ref{lip}). Thus,
\begin{equation}
\label{lipl2}  
\|F(v_{1}(x,s), x)-F(v_{2}(x,s), x)\|_{L^{2}({\mathbb R})}\leq l
\|v_{1}(x,s)-v_{2}(x,s)\|_{L^{2}({\mathbb R})},
\end{equation}
such that
$$
\|e^{(t-s)\{-p^{6}+ibp+a\}}\sqrt{2\pi}\widehat{G}(p)
[\widehat{f}_{v_{1}}(p,s)-\widehat{f}_{v_{2}}(p,s)]\|_{L^{2}({\mathbb R})}\leq
$$
\begin{equation}
\label{eatgvl2}  
e^{aT}l\|G(x)\|_{L^{1}({\mathbb R})}\|v_{1}(x,s)-v_{2}(x,s)\|_{L^{2}({\mathbb R})}.   
\end{equation}
Using (\ref{int0tetsf12}) along with (\ref{eatgvl2}), we derive
$$
\|\widehat{u_{1}}(p,t)-\widehat{u_{2}}(p,t)\|_{L^{2}({\mathbb R})}\leq
$$
$$
e^{aT}l\|G(x)\|_{L^{1}({\mathbb R})}\int_{0}^{T}
\|v_{1}(x,s)-v_{2}(x,s)\|_{L^{2}({\mathbb R})}ds.
$$
By means of the Schwarz inequality 
\begin{equation}
\label{sch2}  
\int_{0}^{T}\|v_{1}(x,s)-v_{2}(x,s)\|_{L^{2}({\mathbb R})}ds\leq
\sqrt{\int_{0}^{T}\|v_{1}(x,s)-v_{2}(x,s)\|_{L^{2}({\mathbb R})}^{2}ds}\sqrt{T}.
\end{equation}
Therefore,
$$  
\|\widehat{u_{1}}(p,t)-\widehat{u_{2}}(p,t)\|_{L^{2}({\mathbb R})}\leq
$$
\begin{equation}
\label{u1hu2hl2}  
e^{aT}l\sqrt{T}\|G(x)\|_{L^{1}({\mathbb R})}\|v_{1}(x,t)-v_{2}(x,t)\|_
{L^{2}({\mathbb R}\times [0, T])}.
\end{equation}
Clearly,
$$
\|u_{1}(x,t)-u_{2}(x,t)\|_{L^{2}({\mathbb R}\times [0, T])}^{2}=\int_{0}^{T}
\|\widehat{u_{1}}(p,t)-\widehat{u_{2}}(p,t)\|_{L^{2}({\mathbb R})}^{2}dt\leq
$$
\begin{equation}
\label{u1u2l2v1v2}
e^{2aT}l^{2}T^{2}\|G(x)\|_{L^{1}({\mathbb R})}^{2}
\|v_{1}(x,t)-v_{2}(x,t)\|_{L^{2}({\mathbb R}\times [0, T])}^{2}.
\end{equation}  
Let us use (\ref{u1u2hint0t}). Hence,
$$
p^{6}[\widehat{u_{1}}(p,t)-\widehat{u_{2}}(p,t)]=
\int_{0}^{t}
e^{(t-s)\{-p^{6}+ibp+a\}}\sqrt{2\pi}p^{6}\widehat{G}(p)
[\widehat{f}_{v_{1}}(p,s)-\widehat{f}_{v_{2}}(p,s)]ds.
$$
Obviously, we have the estimate from above on the norm 
$$
\|p^{6}[\widehat{u_{1}}(p,t)-\widehat{u_{2}}(p,t)]\|_{L^{2}({\mathbb R})}\leq
$$
\begin{equation}
\label{int0tetsf122}  
\int_{0}^{t}\|
e^{(t-s)\{-p^{6}+ibp+a\}}\sqrt{2\pi}p^{6}\widehat{G}(p)
[\widehat{f}_{v_{1}}(p,s)-\widehat{f}_{v_{2}}(p,s)]\|_{L^{2}({\mathbb R})}ds.
\end{equation}
We recall inequality (\ref{fub1}) . Thus,
$$
\|e^{(t-s)\{-p^{6}+ibp+a\}}\sqrt{2\pi}p^{6}\widehat{G}(p)
[\widehat{f}_{v_{1}}(p,s)-\widehat{f}_{v_{2}}(p,s)]\|_{L^{2}({\mathbb R})}^{2}=
$$
$$
2\pi \int_{-\infty}^{\infty}e^{-2(t-s)p^{6}}e^{2(t-s)a}|p^{6}\widehat{G}(p)|^{2}
|\widehat{f}_{v_{1}}(p,s)-\widehat{f}_{v_{2}}(p,s)|^{2}dp\leq 
$$
$$
2\pi e^{2aT}\|p^{6}\widehat{G}(p)\|_{L^{\infty}({\mathbb R})}^{2}\int_{-\infty}^{\infty}
|\widehat{f}_{v_{1}}(p,s)-\widehat{f}_{v_{2}}(p,s)|^{2}dp\leq 
$$
$$
e^{2aT}\Big\|\frac{d^{6}G}{dx^{6}}\Big\|_{L^{1}({\mathbb R})}^{2}
\|F(v_{1}(x,s), x)-F(v_{2}(x,s), x)\|_{L^{2}({\mathbb R})}^{2}.
$$
By virtue of (\ref{lipl2}), we arrive at
$$
\|e^{(t-s)\{-p^{6}+ibp+a\}}\sqrt{2\pi}p^{6}\widehat{G}(p)
[\widehat{f}_{v_{1}}(p,s)-\widehat{f}_{v_{2}}(p,s)]\|_{L^{2}({\mathbb R})}\leq 
$$
\begin{equation}
\label{etsp2gpeat}
e^{aT}l\Big\|\frac{d^{6}G}{dx^{6}}\Big\|_{L^{1}({\mathbb R})}\|v_{1}(x,s)-v_{2}(x,s)\|_
{L^{2}({\mathbb R})}.
\end{equation}
Combining bounds  (\ref{int0tetsf122}), (\ref{etsp2gpeat}) and  (\ref{sch2}) gives us
$$
\|p^{6}[\widehat{u_{1}}(p,t)-\widehat{u_{2}}(p,t)]\|_{L^{2}({\mathbb R})}\leq
$$
\begin{equation}
\label{p2u1hu2hpt}  
e^{aT}\sqrt{T}l\Big\|\frac{d^{6}G}{dx^{6}}\Big\|_{L^{1}({\mathbb R})}
\|v_{1}(x,t)-v_{2}(x,t)\|_{L^{2}({\mathbb R}\times [0, T])}.
\end{equation}
This means that
$$
\Big\|\frac{\partial^{6}}{\partial x^{6}}[u_{1}(x,t)-u_{2}(x,t)]\Big\|_
{L^{2}({\mathbb R}\times [0, T])}^{2}=\int_{0}^{T}
\|p^{6}[\widehat{u_{1}}(p,t)-\widehat{u_{2}}(p,t)]\|_{L^{2}({\mathbb R})}^{2}dt\leq     $$
\begin{equation}
\label{u1u2l2v1v22}
e^{2aT}l^{2}T^{2}\Big\|\frac{d^{6}G}{dx^{6}}\Big\|_{L^{1}({\mathbb R})}^{2}
\|v_{1}(x,t)-v_{2}(x,t)\|_{L^{2}({\mathbb R}\times [0, T])}^{2}.
\end{equation}  
Let us recall (\ref{u1u2hint0t}). Hence,
$$
\frac{\partial}{\partial t}[\widehat{u_{1}}(p,t)-\widehat{u_{2}}(p,t)]=
$$
$$
\{-p^{6}+ibp+a\}[\widehat{u_{1}}(p,t)-\widehat{u_{2}}(p,t)]+
\sqrt{2\pi}\widehat{G}(p)[\widehat{f}_{v_{1}}(p,t)-\widehat{f}_{v_{2}}(p,t)].
$$
Evidently,
$$
\Big\|\frac{\partial}{\partial t}[\widehat{u_{1}}(p,t)-\widehat{u_{2}}(p,t)]
\Big\|_{L^{2}({\mathbb R})}\leq a\|\widehat{u_{1}}(p,t)-\widehat{u_{2}}(p,t)\|_
{L^{2}({\mathbb R})}+
$$
$$      
|b|\|p[\widehat{u_{1}}(p,t)-\widehat{u_{2}}(p,t)]\|_
{L^{2}({\mathbb R})}+
\|p^{6}[\widehat{u_{1}}(p,t)-\widehat{u_{2}}(p,t)]\|_{L^{2}({\mathbb R})}+    
$$
\begin{equation}
\label{ddtu1hptu2hpt}
\sqrt{2\pi}\|\widehat{G}(p)[\widehat{f}_{v_{1}}(p,t)-\widehat{f}_{v_{2}}(p,t)]\|_
{L^{2}({\mathbb R})}.
\end{equation}  
By means of  (\ref{u1hu2hl2}), the first term in the right side of
(\ref{ddtu1hptu2hpt}) can be easily bounded from above by
\begin{equation}
\label{au1hu2hl2}  
aqe^{aT}\sqrt{T}l\|v_{1}(x,t)-v_{2}(x,t)\|_{L^{2}({\mathbb R}\times [0, T])},
\end{equation}
where $q$ is defined in (\ref{g}).
Let us estimate the norm as
$$
\|p[\widehat{u_{1}}(p,t)-\widehat{u_{2}}(p,t)]\|_{L^{2}({\mathbb R})}^{2}=
$$
$$
\int_{|p|\leq 1}p^{2}|\widehat{u_{1}}(p,t)-\widehat{u_{2}}(p,t)|^{2}dp+
\int_{|p|>1}p^{2}|\widehat{u_{1}}(p,t)-\widehat{u_{2}}(p,t)|^{2}dp\leq 
$$
$$  
\|\widehat{u_{1}}(p,t)-\widehat{u_{2}}(p,t)\|_{L^{2}({\mathbb R})}^{2}+
\|p^{6}[\widehat{u_{1}}(p,t)-\widehat{u_{2}}(p,t)]\|_{L^{2}({\mathbb R})}^{2}.
$$
By virtue of inequalities (\ref{u1hu2hl2}) and  (\ref{p2u1hu2hpt}),
the second term in the right side of (\ref{ddtu1hptu2hpt}) can be bounded from above by
\begin{equation}
\label{bu1hu2hl2}  
|b|qe^{aT}\sqrt{T}l\|v_{1}(x,t)-v_{2}(x,t)\|_{L^{2}({\mathbb R}\times [0, T])}.
\end{equation}
According to (\ref{p2u1hu2hpt}), the third term in the
right side of (\ref{ddtu1hptu2hpt}) can be estimated from above by
\begin{equation}
\label{bu1hu2hl3}  
qe^{aT}\sqrt{T}l\|v_{1}(x,t)-v_{2}(x,t)\|_{L^{2}({\mathbb R}\times [0, T])}.
\end{equation}
Using (\ref{fub}) and (\ref{lipl2}), we derive that
$$
2\pi\int_{-\infty}^{\infty}|\widehat{G}(p)|^{2}
|\widehat{f}_{v_{1}}(p,t)-\widehat{f}_{v_{2}}(p,t)|^{2}dp\leq
$$
$$
2\pi
\|\widehat{G}(p)\|_{L^{\infty}({\mathbb R})}^{2}\int_{-\infty}^{\infty}
|\widehat{f}_{v_{1}}(p,t)-\widehat{f}_{v_{2}}(p,t)|^{2}dp\leq
$$
$$
\|G(x)\|_{L^{1}({\mathbb R})}^{2}\|F(v_{1}(x,t), x)-F(v_{2}(x,t), x)\|_
{L^{2}({\mathbb R})}^{2}\leq
$$
$$
\|G(x)\|_{L^{1}({\mathbb R})}^{2}l^{2}\|v_{1}(x,t)-v_{2}(x,t)\|_{L^{2}({\mathbb R})}^{2}.  
$$
Then the fourth term in the right side of (\ref{ddtu1hptu2hpt}) can be bounded
from above by
\begin{equation}
\label{bu1hu2hl4}  
ql\|v_{1}(x,t)-v_{2}(x,t)\|_{L^{2}({\mathbb R})}.
\end{equation}
We combine (\ref{au1hu2hl2}), (\ref{bu1hu2hl2}), (\ref{bu1hu2hl3}) and
(\ref{bu1hu2hl4}) and arrive at
$$
\Big\|\frac{\partial}{\partial t}[\widehat{u_{1}}(p,t)-\widehat{u_{2}}(p,t)]
\Big\|_{L^{2}({\mathbb R})}\leq
$$
$$
qe^{aT}\sqrt{T}l\{a+|b|+1\}\|v_{1}(x,t)-v_{2}(x,t)\|_
{L^{2}({\mathbb R}\times [0, T])}+
ql\|v_{1}(x,t)-v_{2}(x,t)\|_{L^{2}({\mathbb R})}.
$$
Therefore,
$$
\Big\|\frac{\partial}{\partial t}(u_{1}(x,t)-u_{2}(x,t))\Big\|_
{L^{2}({\mathbb R}\times [0, T])}^{2}=\int_{0}^{T}
\Big\|\frac{\partial}{\partial t}[\widehat{u_{1}}(p,t)-\widehat{u_{2}}(p,t)]
\Big\|_{L^{2}({\mathbb R})}^{2}dt\leq           
$$
\begin{equation}
\label{ddtu1u2l2}  
2q^{2}l^{2}[e^{2aT}T^{2}\{a+|b|+1\}^{2}+1]
\|v_{1}(x,t)-v_{2}(x,t)\|_{L^{2}({\mathbb R}\times [0, T])}^{2}.
\end{equation}
Let us recall the definition of the norm (\ref{142n}) . By means of estimates
(\ref{u1u2l2v1v2}), (\ref{u1u2l2v1v22}) and (\ref{ddtu1u2l2}), we obtain
$$
\|u_{1}-u_{2}\|_{W^{1, (6, 2)}({\mathbb R}\times [0, T])}\leq
$$
\begin{equation}
\label{contr}
ql\sqrt{T^{2}e^{2aT}(1+2[a+|b|+1]^{2})+2}
\|v_{1}-v_{2}\|_{W^{1, (6, 2)}({\mathbb R}\times [0, T])}.
\end{equation}
The constant in the right side of (\ref{contr}) is less than one due to
condition (\ref{qlt}). Then under the given assumptions  equation
(\ref{aux}) defines the map
$$
t_{a, b}: W^{1, (6, 2)}({\mathbb R}\times [0, T])\to W^{1, (6, 2)}
({\mathbb R}\times [0, T]),
$$
which is a strict contraction. Its unique fixed
point $w(x,t)$ is the only solution of
problem (\ref{h}), (\ref{ic}) in
$W^{1, (6, 2)}({\mathbb R}\times [0, T])$. \hfill\lanbox

\bigskip

\noindent
{\it Proof of Corollary 1.4.} The statement of the Corollary holds true, which
comes from the fact that the constant in the right side of estimate
(\ref{contr}) does not depend on the initial condition (\ref{ic}) (see e.g. ~\cite{EOE20}). This means that
problem (\ref{h}), (\ref{ic}) has a unique solution
$w(x,t)\in W^{1, (6, 2)}({\mathbb R}\times {\mathbb R}^{+})$. Let us suppose
that $w(x,t)\equiv 0$ for $x\in {\mathbb R}$ and $t\in {\mathbb R}^{+}$.
This will contradict to our assumption that
$\hbox{supp}\widehat{F(0, x)}\cap \hbox{supp}\widehat{G}$ is a set of
nonzero Lebesgue measure on the real line. \hfill\lanbox

\bigskip


\centerline{\bf 3. Acknowledgement}

\bigskip

\noindent
V.V. is grateful to Israel Michael Sigal for the partial support by the
NSERC grant NA 7901.

\end{document}